# The Definition and Numerical Method of Final Value Problem and Arbitrary Value Problem

Shixiong Wang[1*], Jianhua He[1], Chen Wang[2], Xitong Li[1]

[1]School of Electronics and Information, Northwestern Polytechnical University, Xi'an 710129, China
[2]Department of Computer Science, University College London, London WC1E 6BT, UK

* Corresponding author, e-mail: wsx@mail.nwpu.edu.cn

**ABSTRACT**: Many Engineering Problems could be mathematically described by Final Value Problem, which is the inverse problem of Initial Value Problem. Accordingly, the paper studies the final value problem in the field of ODE problems and analyses the differences and relations between initial and final value problems. The more general new concept of the endpoints-value problem which could describe both initial and final problems is proposed. Further, we extend the concept into inner-interval value problem and arbitrary value problem and point out that both endpoints-value problem and inner-interval value problem are special forms of arbitrary value problem. Particularly, the existence and uniqueness of the solutions of final value problem and inner-interval value problem of first order ordinary differential equation are proved for discrete problems. The numerical calculation formulas of the problems are derived, and for each algorithm, we propose the convergence and stability conditions of them. Furthermore, multivariate and high-order final value problems are further studied, and the condition of fixed delay is also discussed in this paper. At last, the effectiveness of the considered methods is validated by numerical experiment.

**KEYWORDS**: differential equation, inverse problem, endpoints-value problem, inner-interval value problem, numerical method

## 1 Introduction

Numerical analysis of differential equations is of great importance in both computational mathematics and practical engineering. A typical problem pertaining to this area is initial value problem, which aims to solve the differential equation by using the initial value of the system. Many previous studies focusing on both numerical calculation and theoretical analysis have established various mature methods in the area. After the wide usage of the digital computer, solving initial value problem with advanced algorithms become much easier, while the problem itself also becomes an important part of Computer Science. Applications of initial value problem in practical areas are broad, ranging from electronic analysis to aerospace engineering, and it has made significant contribution to modern engineering.

A typical application of initial value problem is the inference of a causal system. Often, it is necessary to predict the future condition from the current value of the system. With the utilization of the concept and numerical calculation method of initial value problem, it is possible to deduce any causal system with a digital computer. However, if we consider the inverse problem of causal system deduction [1, 2], it will be found that existed concepts and methods are not able to analyse or solve it. That means the problem about the derivation of the previous status of the system from the current situation still remains unsolved. The problem stated above is actually not a kind of initial value problem (further discussion on this point will be detailed in Remark 1). Thus, it is necessary to find a new concept and the numerical method of it to solve the problem.

The similar problem could be found in some practical areas such as Aviation Fire Control. In this area, a typical problem is how to obtain the initial weapon control conditions (such as aircraft position, posture and weapon delivery time) by inference from ideal attacking effects. In other fields such as control engineering, there also exists the problem of deducing a former situation of the system from the current condition. This paper aims to discuss the description and numerical method of these engineering problems. Rather than only focusing on describing the practical problems, the paper generalizes the problems and proposes the mathematical concept of final value problem. Subsequently, many similar but distinct mathematical concepts which can define more general numerical problems are proposed. The paper then discusses the properties of the new concepts and deduce the numerical calculation methods of them. Properties of the methods such like stability and convergence are also derived.

In summary, this paper makes two novel contributions. Firstly, the concepts proposed in the paper extend the area of numerical method. With our new concepts, it will be possible to analyse more engineering problems which were once difficult to analyse. Secondly, the paper derives numerical methods and analyses the properties of them. These numerical methods could be directly implemented in other areas, while the analyses of the properties of the new concepts are also pioneering works for the new concepts.

## 2 Arbitrary Value Problem of Ordinary Differential Equations
### 2.1 Endpoints-value Problems and its Extension

Mathematically, the problems stated before pertains to the final value problem of ordinary differential equations. The concept of final value problem of first order differential equations can be given as:

$$\begin{cases} \dfrac{dy}{dx} = f(x,y), a \leq x < b \\ y(b) = y_E \end{cases} \quad (1)$$

where $f(x,y)$ is a real-valued function and $y_E$ is the given final value. The variable of the solution to the problem should be on the interval $[a,b]$ and the solution should satisfy the given conditions. The problem in higher dimensions is similar to the definition.

It is noticeable that final value problem is actually not a special case of initial value problem. From the perspective of Physics, final value problem could present the sense of inverse time, for it could deduce the parameters of a previous condition. Mathematically, the relationship between initial value problem and final value problem is neither like the relationship between Taylor formula and Maclaurin formula, nor similar to the relationship between Trapezoidal formula and second-order Runge-Kutta formula. In fact, final value problem is the inverse problem of initial value problem, although the forms of them are similar. It could be reflected from equations (27), (30), (36), (56) and Remark 1 in the following passages.

Obviously, initial value problem aims to solve the "forward-solution" problem using the initial value of the function, derivative or other conditions that can lead to a solution, while final value problem aims to solve the "backward-solution" problem. From this perspective, we extend the concepts and present a new concept called endpoints-value problem, which can describe both initial and final value problems. In fact, it can be found that both initial value problem and final value problem are special cases of endpoints-value problem. The definition of endpoints-value problem can be given as:

**Definition 1** The endpoint-value problem is an ordinary differential equation problem together with the value on the endpoint(s) (either initial point or final point, or both) of the interval, derivative(s) in any order and other condition(s) that can lead to a solution. The variable of the solution to this problem can be either forward or backward in direction. Here, the forward direction means the incremental direction of the independent variable, while the backward direction means the decreasing direction of it. Noticeably, in some high-dimensional cases, it is possible for some variables to be forward while other variables to be backward.

Particularly, the endpoints-value problem of first order ordinary differential equation can be mathematically described as below.
Known:

$$\begin{cases} \dfrac{dy}{dx} = f(x,y), a \leq x \leq b \\ y(a) = y_0 \ \text{ or } \ y(b) = y_E \end{cases} \quad (2)$$

where $f(x,y)$ is a real-valued function and both $y_0$ and $y_E$ are the given endpoints values. The variable of the solution to the problem should be on the interval $[a,b]$ and the solution should satisfy the given conditions.

### 2.2 Definition of Arbitrary Value Problem

Similar to the definition of endpoints-value problem, it is possible to define a differential equation problem in the interval (without endpoints). The concept of inner-interval is then defined as following.

**Definition 2** The inner-interval value problem is an ordinary differential equation together with value(s) inner the interval, derivative(s) and other conditions which can lead to a solution. The variable of solution to the problem should be on the interval.

Specially, the inner-interval value problem of first order ordinary differential equation could be mathematically described as below.
Known:

$$\begin{cases} \dfrac{dy}{dx} = f(x,y), a \leq x \leq c \\ y(b) = y_I, a < b < c \end{cases} \quad (3)$$

where $f(x,y)$ is a real-valued function and $y_I$ is the given inner-interval value. The variable of solution to the problem should be on the interval $[a,c]$ and the solution should satisfy the given conditions.

Comparing endpoints-value problem with inner-interval value problem, it could be found that the forms of them are somehow similar. That means it is possible to find a more general concept which could both describe endpoints-value problem and inner-interval value problem. The concept of arbitrary value problem, which is the most general problem of the problems stated above, could be given as:

**Definition 3** The arbitrary value problem is an ordinary differential equation problem together with value(s) on any point of the closed interval, derivative(s) and other condition(s) which can lead to a solution. The variable of solution to the problem should be on the interval.

Specially, the arbitrary value problem of first order ordinary differential equation can be mathematically described as below.
Known:

$$\begin{cases} \dfrac{dy}{dx} = f(x,y), a \leq x \leq c \\ y(b) = y_I, a \leq b \leq c \end{cases} \tag{4}$$

where $f(x,y)$ is a real-valued function and $y_I$ is the given value at any point of the interval. The variable of solution to the problem should be on the interval $[a,c]$ and the solution should satisfy the given conditions.

As stated before, we now have the concept of arbitrary value problem. The following theorem describes the existence of arbitrary value problem.

**Theorem 1** *For the first order differential equation, if the solution of both initial value problem and final value problem exist, then the solution of arbitrary value problem of the differential equation will exist. Specifically, if the value of the given condition could be on the endpoints, then existence condition of the solution of arbitrary value problem would be weakened to the existence of the solution of the corresponding endpoints-value problem.*

Obviously, the theorem above can be established because the arbitrary value problem can be converted into the combination of initial value problem and final value problem, while both of them can be solved respectively. Specially, the function $f(x,y)$ of inner-interval value problem is not necessarily strictly continuous or smooth on the interval (Fig. 1). This is because initial value problem and final value problem are independent of each other and therefore it is possible to solve the inner-interval value problem of a piecewise function if the value on the discontinuity point(s) could be obtained.

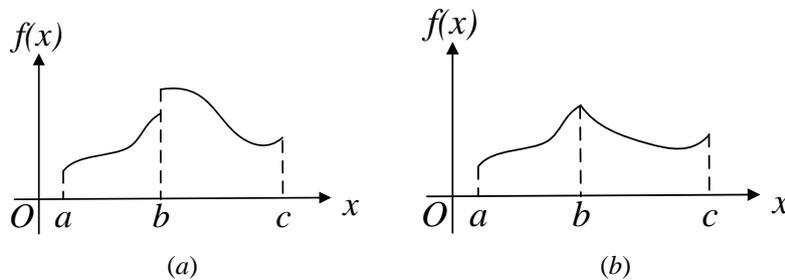

(a)    (b)

**Fig. 1** The function of arbitrary value problem is not necessarily strictly continuous or smooth. From the figure, it can be found that it is possible to solve the inner-interval value problem of both functions shown in the figure if the value on *b* point can be obtained

The theorem of the existence of the solution of initial value problem and its numerical method have already been established by previous studies [3], hence in the following passages our aim is to discuss the existence of the solution of finial value problem and find the numerical calculation method(s) of it.

**3 The Existence of the Solution of Final Value Problem**

**Definition 4** (given in [3]) A function $f(x,y)$ is said to satisfy the Lipschitz condition in the variable $y$ if a constant $L > 0$ exists with

$$|f(x,y_1) - f(x,y_2)| \leq L|y_1 - y_2| \tag{5}$$

whenever $(x,y_1), (x,y_2) \in D$. The constant $L$ is called a Lipschitz constant for $f(x,y)$.

**Theorem 2** *For the final value problem of first order ordinary differential equation (1):*

*Suppose* $D = \{(x,y) | a \leq x \leq b, -\infty < y < +\infty\}$, *if the function* $f(x,y)$ *could satisfy the following conditions:*

(a) $f(x,y)$ *is continuous on the set* $D$.
(b) $f(x,y)$ *satisfies the Lipschitz condition in the variable* $y$ *on the set* $D$.

*Then the unique solution of final value problem exists. In addition, the solution is continuous and differential. It could also be found that the solution is dependent on the function and the final value given as the condition.*

*Proof:* The theorem can be proved by Picard's Existence and Uniqueness Theorem [3]. □

Intuitively, the theorem is a natural deduction from the existence and uniqueness of initial value problem. Comparing the form of initial value problem with the form of final value problem, it could be found that they are essentially similar. Using the concept of right derivative it is possible to define and solve initial value problem, while in the same way, it is possible to define and solve final value problem using the concept of left derivative. The idea will be discussed further in the following passages.

**Remark 1** It should be noted that we could NOT consider the final value problem (1) as one kind of initial value problem just by changing the independent variable *x* to *-x*.

If we let *-t = x* in (1), we can then get

$$\frac{dy}{dt} = -f(-t, y), -b < t \leq -a \quad (6)$$

but how can we decide the initial value of initial value problem? Namely, *y(-b) = what*? Generally, NO ANY effective methods could be taken to get that!

In fact, for a practical or engineering problem, if *t* stands for time and the point *b* presents the current time point, then, we can get the system state of time point *b* definitely does not mean that we can also get the system state of the time point *-b*, thus no any efficient steps can be taken to convert a final value problem to an equivalent initial value problem, it is not possible to do that by only replacing *t* to *-t* at least. Thus the final value problem is essentially different from the initial value problem.

**4 Numerical Solution Methods of Final Value Problem**

According to the definition of derivative

$$\frac{dy}{dx}\Big|_{x=x_0} = \lim_{\Delta x \to 0} \frac{y(x_0 + \Delta x) - y(x_0)}{\Delta x} \quad (7a)$$
$$= \lim_{\Delta x \to 0} \frac{y(x_0) - y(x_0 - \Delta x)}{\Delta x} \quad (7b)$$
$$= \lim_{\Delta x \to 0} \frac{1}{2}\frac{y(x_0 + \Delta x) - y(x_0 - \Delta x)}{\Delta x} \quad (7c)$$

(7)

obviously, the sufficient and necessary condition of the existence of derivative is the existence and equality of both left derivative and right derivative.

In order to deduce the numerical solution method, two forms of Taylor formula could be taken into account.

$$\begin{cases} y(x_{n+1}) = y(x_n) + hy^{(1)}(x_n) + \cdots + \frac{1}{n!}h^n y^{(n)}(x_n) + o(h^n) & (8a) \\ y(x_n) = y(x_{n+1}) - hy^{(1)}(x_{n+1}) + \cdots + \frac{1}{n!}h^n y^{(n)}(x_{n+1}) + o(h^n) & (8b) \end{cases}$$

(8)

Equation (7*a*) is actually the fundament of the recurrence formula of initial value problem. Replacing $\Delta x$ with a sufficient precise step size *h* is the mean idea of the recurrence formula. Thus it is possible to use the discrete difference operation to replace continuous differential operation. Moreover, equation (7*c*) and (8*a*) is the fundament of the numerical method with higher precision. The idea of the replacement is the origin of most classical numerical methods of initial value problem such as Euler's method or Trapezoidal method. Naturally, the idea that numerical method of final value problem could also be deduced from the same approach is then proposed. In fact, the numerical solution method of final value problem in this paper does come from equation (7*b*) and (8*a*). The relevant deduction will be illustrated in the following passages.

In order to accent the backward solution of the problem, the recurrence formula will be wrote in the form of $y_n = T(y_{n+1})$, where $T(\cdot)$ is the recursion operator. Here $y_n$ stands for the approximation of $n_{th}$ order derivative while $y_{n+1}$ indicates the $(n-1)_{th}$ result from the final value, showing that the order of the steps is backward (the first step is from the final point).

## 4.1 The Global Convergence Theorem of Numerical Solution Method of Final Value Problem

**Definition 5** The definition of the general form of one-step recurrence formula is shown as below
$$y_n = y_{n+1} - h\varphi(x_{n+1}, y_{n+1}, h) \tag{9}$$

where $\varphi(x, y, h)$ is called **Decrement Function**.

**Theorem 3** *Suppose the one-step method has order p and the given final value is precise. If decrement function $\varphi(x, y, h)$ satisfies a Lipschitz condition in variable y:*
$$|\varphi(x, y_1, h) - \varphi(x, y_2, h)| \leq L|y_1 - y_2| \tag{10}$$

*then the global truncation error will satisfy:*
$$|e_n| = o(h^p) \tag{11}$$

*where L is the Lipschitz constant.*

*Proof :* Assume that $y_{n+1} = y(x_{n+1})$ in the numerical method of final value problem. $\bar{y}_n$ is then used to indicate the result. Thus we have the following equation:
$$\bar{y}_n = y(x_{n+1}) - h\varphi(x_{n+1}, y(x_{n+1}), h) \tag{12}$$

According to the definition, $|y(x_n) - \bar{y}_n|$ is the local truncation error. Since the method has order $p$, there exists a positive constant $C$ to satisfy the following inequality.
$$|y(x_n) - \bar{y}_n| \leq Ch^{p+1} \tag{13}$$

Furthermore, according to equation (9) and (12), following inequality could be established.
$$|\bar{y}_n - y_n| \leq |y(x_{n+1}) - y_{n+1}| + h|\varphi(x_{n+1}, y(x_{n+1}), h) - \varphi(x_{n+1}, y_{n+1}, h)| \tag{14}$$

Now using the Lipschitz Condition, $|\bar{y}_n - y_n|$ could satisfy the following inequality.
$$|\bar{y}_n - y_n| \leq (1 + hL)|y(x_{n+1}) - y_{n+1}| \tag{15}$$

Thus we can get the inequality:
$$|y(x_n) - y_n| \leq |\bar{y}_n - y_n| + |y(x_n) - \bar{y}_n| \leq (1 + hL)|y(x_{n+1}) - y_{n+1}| + Ch^{p+1} \tag{16}$$

That means, for the global truncation error $|e_n| = |y(x_n) - y_n|$, the following recurrence formula could be established:
$$|e_n| \leq (1 + hL)|e_{n+1}| + Ch^{p+1} \tag{17}$$

Then the inequality of global truncation error should be:
$$|e_n| \leq (1 + hL)^n|e_0| + \frac{Ch^p}{L}[(1 + hL)^n - 1] \tag{18}$$

where $e_0$ is the initial error. It can also be noted as $e_E$, although $e_0$ is more suitable for understanding.
Noticeably, if $x_E - x_n = nh \leq T$, then:
$$(1 + hL)^n \leq (e^{hL})^n \leq e^{TL} \tag{19}$$

Thus the ultimate estimation formula could be derived:
$$|e_n| \leq |e_0|e^{TL} + \frac{Ch^p}{L}(e^{TL} - 1) \tag{20}$$

Consequently, the conclusion is said to be correct if $|e_0| = 0$ (that is, the given value is precise). □

In other words, the convergence of one-step method dependents on whether the decrement functions could satisfy a Lipschitz Condition.

## 5 Basic Numerical Methods of Final Value Problem
### 5.1 Method in Euler's Form and its Convergence and Stability Condition

Similar to the recurrence formula of initial value problem introduced in [3], the recurrence formula of Euler's form in final value problem is established below.
Explicit Euler's form difference equation

$$y_n = y_{n+1} - hf(x_{n+1}, y_{n+1})$$
(21)

Implicit Euler's form difference equation
$$y_n = y_{n+1} - hf(x_n, y_n)$$
(22)

For the Implicit case, usually the explicit formula will first be used to make a prediction. Then the implicit formula will be used to calculate a more precise value. This is called Predictor-Corrector method, the formula of it is illustrated below.
$$\begin{cases} y_n^{(0)} = y_{n+1} - hf(x_{n+1}, y_{n+1}) \\ y_n^{(s+1)} = y_{n+1} - hf(x_n, y_n^{(s)}) \end{cases}$$
(23)

where $s = 1,2,3,...$ means the iteration times.

Now the convergence and stability of the solution will be discussed. If we use equation (23) to subtract equation (22), with necessary deduction, we can get the following equation.
$$|y_n^{(s+1)} - y_n| = h|f(x_n, y_n^{(s)}) - f(x_n, y_n)|$$
(24)

According to Theorem 2 (Lipschitz condition)
$$|y_n^{(s+1)} - y_n| = h|f(x_n, y_n^{(s)}) - f(x_n, y_n)|$$
$$\leq hL|y_n^{(s)} - y_n|$$
(25)

Then from the inequality we can see that the convergence condition of the method in the form of Predictor-Corrector (Equation 23) is
$$hL < 1$$
(26)

where $L$ is the Lipschitz Constant, and $h$ is the step size.

For the stability, the common model of stability analysis is used.
$$\frac{dy}{dx} = \lambda y, \lambda \in R^+$$
(27)

Notice that the condition of $\lambda$ is different from the condition in initial value problem. In initial value problem the condition should be $\lambda \in R^-$.

Take the explicit Euler's form method as example, the stability condition of it can be proved below. Using the common model, equation (28) can be established.
$$y_i = y_{i+1} - hf(x_{i+1}, y_{i+1}) = y_{i+1} - h\lambda y_{i+1} = (1 - h\lambda)y_{i+1}$$
(28)

Take the errors into account when calculating $y_{i+1}$, it can be represented as:
$$\rho_{i+1} = \bar{y}_{i+1} - y_{i+1}$$
where $\bar{y}_{i+1}$ is the approximation of $y_{i+1}$ ($y_{i+1}$ is the theoretical value of $\bar{y}_{i+1}$ ).

Then we can obtain the equation:
$$|\bar{y}_i - y_i| = |1 - h\lambda||\bar{y}_{i+1} - y_{i+1}|$$
(29)

where $\bar{y}_i$ is the approximation of $y_i$.

Consequently, from equation (29) the (necessary) stability condition of explicit Euler's form method can be obtained. The condition is:
$$|1 - h\lambda| < 1$$
(30)

**5.2 Method in Trapezoidal Form and its Convergence and Stability Condition**

In the same way of Euler's method, the recurrence formula of method in Trapezoid form could be established.
$$y_n = y_{n+1} - \frac{h}{2}[f(x_n, y_n) + f(x_{n+1}, y_{n+1})]$$
(31)

Together with explicit Euler's form method, the following combined Predictor-Corrector method could be established.
$$\begin{cases} y_n^0 = y_{n+1} - hf(x_{n+1}, y_{n+1}) \\ y_n^{(s+1)} = y_{n+1} - \frac{h}{2}[f(x_n, y_n^{(s)}) + f(x_{n+1}, y_{n+1})] \end{cases}$$

(32)

Similar to what has been discussed in the section of Euler's method, if we use equation (32) to subtract equation (31) with some derivation, the following equation can be obtained.

$$|y_n^{(s+1)} - y_n| = \frac{h}{2}|f(x_n, y_n^{(s)}) - f(x_n, y_n)|$$

(33)

According to Theorem 2 (Lipschitz Condition)

$$|y_n^{(s+1)} - y_n| = \frac{h}{2}|f(x_n, y_n^{(s)}) - f(x_n, y_n)|$$
$$\leq \frac{h}{2}L|y_n^{(s)} - y_n|$$

(34)

Then, the convergence condition of the algorithm proposed in equation (32) (called Predictor-Corrector method in trapezoid form) is:

$$\frac{h}{2}L < 1$$

(35)

where $L$ is Lipschitz constant, and $h$ is the step size.

In the same way of Euler's method, the stable condition of Trapezoid method can be obtained.

$$\left|\frac{1 - \frac{\lambda}{2}h}{1 + \frac{\lambda}{2}h}\right| < 1$$

(36)

The inequality always holds because $\lambda > 0$. Thus the trapezoid method (31) of final value problem is unconditionally stable.

### 5.3 Runge-Kutta Method of Final Value Problem

Similar to the method in initial value problem, the formula of Runge-Kutta Method in final value problem is then established.

$$\begin{cases} y_n = y_{n+1} - h\sum_{i=1}^{r} c_i K_i \\ K_1 = f(x_{n+1}, y_{n+1}) \\ K_i = f\left(x_{n+1} - \alpha_i h, y_{n+1} - h\sum_{j=1}^{i-1}\beta_j K_j\right) \end{cases}$$

(37)

Equation (37) is a general form. With different value of $r$, $c$, $\alpha$ and $\beta$ it will have different meaning and properties.

### 5.4 Convergence and Stability Condition of Fourth-order Classical Runge-Kutta Method

The algorithm of Fourth-order Classical Runge-Kutta Method in final value problem will be discussed and deduced below.

It will be extremely complicate if we directly deduce the forth-order Runge-Kutta formula from equation (37). Alternatively, if we consider the integral form of equation (1):

$$y_n = y_{n+1} - \int_{x_n}^{x_{n+1}} f(x, y)dx$$

(38)

The target Runge-Kutta formula has order 4 in precision, thus it will be feasible to approximate the integral $\int_{x_n}^{x_{n+1}} f(x, y)dx$ by using a numerical method with truncation error of $O(h^5)$. The Simpson's rule [3] is suitable for this condition.

In order to obtain the Runge-Kutta formula with the similar form of that in initial value problem, the Simpson's rule in the form of quadratic interpolation is utilized. The equation is shown below:

$$\int_{x_n}^{x_{n+1}} f(x)dx = \frac{h}{6}\left[f(x_n) + 2f\left(\frac{x_n + x_{n+1}}{2}\right) + f(x_{n+1})\right] + \varepsilon$$

(39)

where the cut-off error is $\varepsilon = -\frac{1}{90}f^{(4)}(\xi, y)h^5$, and $\varepsilon \in (x_n, x_{n+1})$.

Consequently, the form of RK4 in final value problem can be illustrated as below.

$$\begin{cases} y_n = y_{n+1} - \dfrac{h}{6}(K_1 + 2K_2 + 2K_3 + K_4) \\ K_1 = f(x_{n+1}, y_{n+1}) \\ K_2 = f\left(x_{n+1} - \dfrac{1}{2}h, y_{n+1} - \dfrac{h}{2}K_1\right) \\ K_3 = f\left(x_{n+1} - \dfrac{1}{2}h, y_{n+1} - \dfrac{h}{2}K_2\right) \\ K_4 = f(x_{n+1} - h, y_{n+1} - h(aK_1 + bK_2 + cK_3)) \end{cases}$$

(40)

where $a$, $b$ and $c$ are undetermined coefficients and subject to the condition $a + b + c = 1$.

To determine the unknown coefficients, the stability analysis equation (Equation 27) is then used:

$$\begin{cases} K_1 = \lambda y_{n+1} \\ K_2 = \lambda\left(1 - \dfrac{1}{2}\lambda h\right) y_{n+1} \\ K_3 = \lambda\left(1 - \dfrac{1}{2}\lambda h + \dfrac{1}{4}\lambda^2 h^2\right) y_{n+1} \\ K_4 = \lambda\left[1 - (a + b + c)\lambda h + \left(\dfrac{b+c}{2}\right)\lambda^2 h^2 - \dfrac{1}{4}c\lambda^3 h^3\right] y_{n+1} \end{cases}$$

(41)

Then we can obtain the following equation:

$$y_n = y_{n+1} - \dfrac{1}{6}\left[6\lambda h - (2 + a + b + c)\lambda^2 h^2 + \dfrac{1}{2}(1 + b + c)\lambda^3 h^3 - \dfrac{1}{4}c\lambda^4 h^4\right] y_{n+1}$$

(42)

Comparing with the corresponding Taylor's formula:

$$y_n = y_{n+1} - \left[\lambda h - \dfrac{1}{2!}\lambda^2 h^2 + \dfrac{1}{3!}\lambda^3 h^3 - \dfrac{1}{4!}\lambda^4 h^4\right] y_{n+1}$$

(43)

The following equations can be obtained:

$$\begin{cases} a + b + c = 1 \\ b + c = 1 \\ c/4 = 1/4 \end{cases}$$

(44)

thus $a = b = 0, c = 1$.

Consequently, the algorithm of Fourth-order Classical Runge-Kutta Method in final value problem can be shown as below:

$$\begin{cases} y_n = y_{n+1} - \dfrac{h}{6}(K_1 + 2K_2 + 2K_3 + K_4) \\ K_1 = f(x_{n+1}, y_{n+1}) \\ K_2 = f\left(x_{n+1} - \dfrac{1}{2}h, y_{n+1} - \dfrac{h}{2}K_1\right) \\ K_3 = f\left(x_{n+1} - \dfrac{1}{2}h, y_{n+1} - \dfrac{h}{2}K_2\right) \\ K_4 = f(x_{n+1} - h, y_{n+1} - hK_3) \end{cases}$$

(45)

It could be found that the formula of RK4 in final value problem is extremely similar to the formula in initial value problem. The only difference between them is some signs that indicate positive and negative.

The decrement function in algorithm (45) is:

$$\varphi(x, y, h) = \dfrac{1}{6}[k_1(x, y, h) + 2k_2(x, y, h) + 2k_3(x, y, h) + k_4(x, y, h)]$$

(46)

where

$$\begin{cases} k_1(x, y, h) = f(x, y) \\ k_2(x, y, h) = f\left(x - \dfrac{1}{2}h, y - \dfrac{1}{2}hk_1(x, y, h)\right) \\ k_3(x, y, h) = f\left(x - \dfrac{1}{2}h, y - \dfrac{1}{2}hk_2(x, y, h)\right) \\ k_4(x, y, h) = f(x - h, y - hk_3(x, y, h)) \end{cases}$$

(47)

According to Theorem 2 (Lipschitz Condition)

$$|k_1(x, y_1, h) - k_1(x, y_2, h)| \le L|y_1 - y_2|$$

$$|k_2(x, y_1, h) - k_2(x, y_2, h)| \leq L \left( |y_1 - y_2| + \frac{1}{2} hL|y_1 - y_2| \right) \tag{48}$$
$$= L(1 + \frac{1}{2}h)|y_1 - y_2| \tag{49}$$

$$|k_3(x, y_1, h) - k_3(x, y_2, h)| \leq L \left( |y_1 - y_2| + \frac{1}{2} hL \left(1 + \frac{1}{2}h\right) |y_1 - y_2| \right)$$
$$= L(1 + \frac{1}{2} hL + (\frac{1}{2} hL)^2)|y_1 - y_2| \tag{50}$$

$$|k_4(x, y_1, h) - k_4(x, y_2, h)| \leq L \left( |y_1 - y_2| + hL \left(1 + \frac{1}{2} hL + \left(\frac{1}{2} hL\right)^2 \right) |y_1 - y_2| \right)$$
$$= L(1 + hL + \frac{1}{2}(hL)^2 + \frac{1}{4}(hL)^3)|y_1 - y_2| \tag{51}$$

In the deduction, the inequality
$$|a - b| \leq |a| + |b| \qquad a, b \in R \tag{52}$$

is used. This is similar to the corresponding deduction in initial value problem. In that process the similar inequality shown below is used.
$$|a + b| \leq |a| + |b| \qquad a, b \in R \tag{53}$$

As derived above, we can obtain the inequality:
$$|\varphi(x, y_1, h) - \varphi(x, y_2, h)| \leq L \left[ 1 + \frac{1}{2} hL + \frac{1}{6} h^2 L^2 + \frac{1}{24} h^3 L^3 \right] |y_1 - y_2| \tag{54}$$

Assuming that $\widetilde{L} = L \left[ 1 + \frac{1}{2} hL + \frac{1}{6} h^2 L^2 + \frac{1}{24} h^3 L^3 \right]$, according to Theorem 3, the algorithm of final value problem is said to be convergent if $\widetilde{L} > 0$. Obviously, the condition could be satisfied by an appropriate $hL$.
The solution of the inequality $\widetilde{L} > 0$ is:
$$hL > -2.7853 \tag{55}$$

The inequality (55) always holds. Consequently, Classical RK4 in final value problem is unconditional convergent.
The stability condition of Classical RK4 for final value problem can later be deduced by analyzing algorithm (45) with the common model of stability analysis (Equation 27).
The condition is:
$$\left| 1 - \lambda h + \frac{1}{2} \lambda^2 h^2 - \frac{1}{6} \lambda^3 h^3 + \frac{1}{24} \lambda^4 h^4 \right| < 1 \tag{56}$$

**5.5 Runge-Kutta Method in Other Orders**
Using the method of deduction stated above, it is easy to derive the algorithm and the stability and convergence condition of Runge-Kutta methods in other orders.

**6 The Numerical Method of Multivariate Final Value Problem and High-order Final Value Problem**
The multivariable differential equations could be given as:
$$\begin{cases} \frac{dy_i}{dx} = f_i(x, y_1, y_2, \cdots, y_i, \cdots, y_n), a \leq x < b \\ y_i(b) = y_{iE} \end{cases} \tag{57}$$

where $i = 1, 2, 3, \cdots, n$ indicates the $i^{th}$ differential equation.
Similar with the derivation of single-variable Initial Value Problem, the numerical method of multivariate Final Value Problem in the form of (57) could be derived with the theorem proposed in this paper. The Fourth-order Runge-Kutta method of high-order differential equation in the form of (57) could be given as:

$$\begin{cases} y_{im} = y_{im+1} - \frac{h}{6}(K_{i1} + 2K_{i2} + 2K_{i3} + K_{i4}) \\ K_{i1} = f_i(x_{m+1}, y_{1m+1}, \cdots, y_{nm+1}) \\ K_{i2} = f_i\left(x_{m+1} - \frac{h}{2}, y_{1m+1} - \frac{h}{2}K_{11}, \cdots, y_{nm+1} - \frac{h}{2}K_{n1}\right) \\ K_{i3} = f_i\left(x_{m+1} - \frac{h}{2}, y_{1m+1} - \frac{h}{2}K_{12}, \cdots, y_{nm+1} - \frac{h}{2}K_{n2}\right) \\ K_{i4} = f_i(x_{m+1} - h, y_{1m+1} - hK_{13}, \cdots, y_{nm+1} - hK_{n3}) \end{cases}$$

(58)

where $i = 1,2,3,\cdots,n$ indicates the $i^{th}$ differential equation and $m$ marks the $m^{th}$ step of iteration.

Similarly, we can also formulate the High-order Final Value Problem. For universality and simplicity, consider only the high-order differential equation in general polynomial form:

$$y^{(n)} + p_{n-1}y^{(n-1)} + \cdots + p_1 y + p_0 = f$$

(59)

where $p_i, (i = 0,1,2,\dots,n-1)$ is the coefficient polynomial function of $y^{(i)}$, and $y^{(i)}$ stands for $d^i y / dx^i$.

It could be argued that it is always possible to transfer the differential equations in the form of (59) to the form of (57) by defining suitable and enough state variables [4, 5]. Consequently, the numerical method in (58) is also suitable for general high-order Final Value Problem. An instance will be illustrated in Example 1.

**Example 1** Example of Converting a High-order Final Value Problem to its Corresponding Multivariate Final Value Problem.

Consider a high-order differential equation in general polynomial form
$$4\dddot{x} + 3\ddot{x} + 2\dot{x} + 1 = 0$$
If we define $x_1 = x; x_2 = \dot{x}; x_3 = \ddot{x}$ ($x_1, x_2$ and $x_3$ is always called State Variable [4, 5]), then we can convert it to the following equivalent form
$$\begin{cases} \dot{x}_3 = (-3x_3 - 2x_2 - 1)/4 \\ \dot{x}_2 = x_3 \\ \dot{x}_1 = x_2 \end{cases}$$

Obviously, the high-order final value problem was converted to its corresponding multivariable final value problem form. Thus if the final value of state variables are given (for example, we can set $x_1 = 0; x_2 = 0; x_3 = 0$), the problem could then be solved by equation (58).

## 7 The Numerical Method of Final Value Problem with Fixed Delay

The delay condition is an important part of initial value problem and has been discussed by many researchers [6, 7]. Likewise, it also exists in Final Value Problem. This section will discuss the numerical method of first order Final value Problem with fixed delay.

The differential equation in this case could be given as:
$$\begin{cases} \frac{dy}{dx} = f(x + T, y), a \leq x < b \\ y(b) = y_E \end{cases}$$

(60)

According to Laplace Transform [4, 5], $f(x + T, y)$ in equation (60) could be transferred into the following form:
$$F(x + T, y)(s) = F(x, y)(s)e^{Ts}$$

(61)

where $F(x, y)(s)$ indicates the Laplace Transform of $f(x, y)$.

Noticeably, in Initial Value Problem, the similar equation should be
$$F(x - T, y)(s) = F(x, y)(s)e^{-Ts}$$

(62)

The difference between (61) and (62) further illustrates the inverse properties between final and initial value problem.

Moreover, according to Taylor's formula, the exponential part of equation (62) could be transferred into:
$$e^{Ts} = 1 + \frac{1}{1!}(Ts) + \frac{1}{2!}(Ts)^2 + \cdots + \frac{1}{n!}(Ts)^n + \varepsilon$$

(63)

where $\varepsilon$ indicates the truncation error.

Consequently, it could be found that the problem proposed with equation (60) could always be transferred into the problem proposed with equation (59) with arbitrarily small error [4, 5]. Therefore the problem could be further transferred into (57), which could be solved by equation (58) proposed in this paper.

## 8 Experimental Result and Analysis
### 8.1 Using RK4 Method in both Initial Value Problem and Final Value Problem to Solve a Classic Problem

Considering following initial value problem and final value problem:

$$\begin{cases} \frac{dy}{dx} = y - \frac{2x}{y} \\ y(0) = 1.0 \end{cases} \text{ and } \begin{cases} \frac{dy}{dx} = y - \frac{2x}{y} \\ y(1) = \sqrt{3} \end{cases}$$

The differential equations of the two problems are same. Suppose $h = 0.1$ and solve the problems with RK4 method in both initial value problem and final value problem (Equation 45) respectively. The results are illustrated in Table 1.

**Table 1** The Numerical Results of The Problem

| Initial Value Problem | | | Final Value Problem | | |
|---|---|---|---|---|---|
| $x_n$ | $y_n$ | $|y_n - y(x_n)|$ | $x_n$ | $y_n$ | $|y_n - y(x_n)|$ |
| 0.0 | 1.000000 | 0.000000 | 1.0 | 1.732051 | 0.000000 |
| 0.1 | 1.095446 | 0.000000 | 0.9 | 1.673320 | 0.000000 |
| 0.2 | 1.183217 | 0.000001 | 0.8 | 1.612451 | 0.000000 |
| 0.3 | 1.264912 | 0.000001 | 0.7 | 1.549193 | 0.000000 |
| 0.4 | 1.341642 | 0.000002 | 0.6 | 1.483239 | 0.000000 |
| 0.5 | 1.414216 | 0.000002 | 0.5 | 1.414213 | 0.000001 |
| 0.6 | 1.483242 | 0.000003 | 0.4 | 1.341640 | 0.000001 |
| 0.7 | 1.549196 | 0.000003 | 0.3 | 1.264910 | 0.000001 |
| 0.8 | 1.612455 | 0.000004 | 0.2 | 1.183215 | 0.000001 |
| 0.9 | 1.673325 | 0.000005 | 0.1 | 1.095444 | 0.000001 |
| 1.0 | 1.732056 | 0.000006 | 0.0 | 0.999999 | 0.000001 |

From the Table 1, it could be found that RK4 in final value problem has sufficient precision in obtaining the solution. The result substantiates the theorem proposed before. The algorithm will be suitable in practical areas, while it also extends the theorem of numerical method mathematically.

## 9 Conclusion

The paper proposes the concept of arbitrary value problem from the question in practical engineering. Related concepts of endpoints-value problem (including initial value problem and final value problem) and inner-interval value problem are also defined. Additionally, different kinds of numerical solution methods in final value problem are deduced. In mathematics, the establishment of these methods has important significance in extending the area of numerical method.

Because of the limitations on space, numerical methods of final value problem in other forms (such as Multi-step method and Richardson Extrapolation [8, 9]) and complex conditions are not studied. Our long-term research plan is to perfect the theory of final value problem and to apply it in more complicate situations, such as the conditions discussed in [6, 7]. We believe that through our further studies the final value problem could be widely applied in both mathematics and practical engineering.

## 10 Acknowledgements

The authors would like to thank Prof. Suzhen Wang, who gave us many useful comments on the research.